\documentclass[letterpaper, 10 pt, conference]{IEEEconf}  
\pdfoutput=1

\IEEEoverridecommandlockouts                              

\usepackage{graphicx}
\usepackage{amsmath} 
\usepackage{amssymb}  
\usepackage{comment}
\let\proof\relax
\let\endproof\relax
\usepackage{xcolor}
\usepackage{graphicx}
 \usepackage{multirow}
\newtheorem{thm}{Theorem}
\newtheorem{lem}{Lemma}
\newtheorem{prop}{Proposition}
\newtheorem{rem}{Remark}
\renewcommand{\QED}{\QEDopen}

\def\proof{{\itshape Proof: }}
\def\endproof{\hspace*{\fill}~\QED\par\endtrivlist\unskip}

\def\e{\varepsilon}
\def\p{\mathfrak{p}}
\def\E{\mathcal{E}}
\def\L{\mathcal{L}}
\DeclareMathOperator{\tr}{tr}
\DeclareMathOperator{\Co}{Co}
\DeclareMathOperator{\diag}{diag}
\allowdisplaybreaks[4]

\title{\LARGE \bf A Multirate Variational Approach to Nonlinear MPC 
}

\author{Yana Lishkova$^{1}$, Mark Cannon$^{1}$ and Sina Ober-Bl\"{o}baum$^{2}$
\thanks{$^{1}$Department of Engineering Science, University of Oxford, Oxford, OX1 3PJ, UK
        {\tt\small \{yana.lishkova,mark.cannon\}@eng.ox.ac.uk}}%
\thanks{$^{2}$Department of Mathematics, Paderborn University, D-33098 Paderborn, Germany
        {\tt\small sinaober@math.uni-paderborn.de}}%
}

\begin{document}

\maketitle
\pagestyle{plain}

\begin{abstract}

A multirate nonlinear model predictive control (NMPC) strategy is proposed for systems with dynamics and control inputs evolving on different timescales. The proposed multirate formulation of the system model and receding horizon optimal control problem allows larger time steps in the prediction horizon compared to single-rate schemes, providing computational savings while ensuring recursive feasibility. A multirate variational model is used with a tube-based successive linearization NMPC strategy. This allows either Jacobian linearization or linearization using quadratic and linear Taylor series approximations of the Lagrangian and generalized forces respectively, providing alternative means for computing linearization error bounds. The two approaches are shown to be equivalent for a specific choice of approximation points and their structure-preserving properties are investigated. Numerical examples are provided to illustrate the multirate approach, its conservation properties and computational savings. 

\end{abstract}

\section{INTRODUCTION}
Model Predictive Control is an 
optimization-based control strategy for
systems with state and input constraints, and model and measurement uncertainty. A key requirement for robust stability is the existence of a feasible trajectory each time the receding horizon optimal control problem (RHOCP) is solved~\cite{smpc}. However, for multirate systems, which exhibit dynamics and/or controls on different timescales, a small discrete step is needed to allow accurate resolution of fast timescales and account for fast changing controls. Small discrete steps increase the computational burden of the RHOCP solved online and reduce the available solution time.

To tackle this problem, non-uniform grid discretization and move-blocking strategies have been developed. The relative merits of these methods in terms of discretization error, computational load, and stability and feasibility guarantees are reviewed in~\cite{cagienard2007move,10in2019mb, OptimMBNew}. 
As an alternative to these approaches we propose a new multirate receding horizon strategy which uses a multirate 
system model to formulate a novel RHOCP formulation and perform the optimization less frequently without loss of fidelity or recursive feasibility. 

Multirate models (e.g.~\cite{Reading9,gunther2016multirate}) allow slow and fast changing states and controls to be treated on separate micro and macro time grids, each with an appropriate time step. In an optimal control setting this can enable controls to vary on different timescales with well understood discretization error behaviour. 
%
%
We consider the  variational method of \cite{Reading9},
%
which is the basis of the multirate Discrete Mechanics and Optimal Control (DMOC) approach~\cite{Reading8}. 
Without decoupling the equations of motion, the method
allows for the resolution of slow and fast changing states on separate macro and micro time grids, reducing the size of the optimization problem
%
and providing a sparse structure for the Jacobian of state constraints~\cite{Reading452,MyASS2020paper}. 
%
%
Another advantage of the variational approach is its preservation of  qualitative model properties from continuous to discrete domains \cite{Reading51}. Both single-rate DMOC \cite{Reading71} and multirate DMOC employ a discrete version of the Lagrange-d’Alembert principle. The resulting discrete-time model inherits conservation properties of the continuous-time model such as Lagrangian symmetries, and accurately represents energy or momenta for exponentially long times~\cite{Reading35}. 
Both DMOC schemes were originally implemented via nonlinear programming (NLP) and as such, single-rate DMOC has previously been used in MPC (e.g.~\cite{Reading785, Reading782}). A naive multirate extension would solve the RHOCP using multirate DMOC at each micro time step. However, depending on the RHOCP size,
the optimisation might not be solvable within this step, despite the savings multirate DMOC provides. 

Instead we propose a multirate receding horizon strategy in which the RHOCP is formulated and solved only at the macro time nodes. This allows larger prediction time steps than single-rate methods, providing computational savings while retaining high-fidelity resolutions of the fast dynamics. Given accuracy and recursive feasibility considerations, we incorporate multirate DMOC within the successive linearisation NMPC scheme of \cite{OrigPaper}. 
This provides stability and recursive feasibility guarantees by bounding the effect of linearization errors using robust tubes. The resulting RHOCP is a convex program, providing significant computational savings.  
%
%
%
To retain the structure-preserving properties of the variational system model, we explore linearization using quadratic and linear approximations of the Lagrangian and generalized forces respectively. 
We show that this approach is equivalent to Jacobian linearization of the discrete variational model for particular choices of approximation points and derive corresponding error bounds. 
To our knowledge this approach has not been used before with variational models.

The paper is organised as follows: Section~\ref{sec:modelling} presents the single-rate and multirate variational model schemes and Section~\ref{sec:linearization} investigates the proposed linearization approach. Section~\ref{sec:mpc} poses the optimal control problem and summarizes the tube-based NMPC strategy. Section~\ref{sec:multirate_mpc} brings these topics together in a novel multirate receding horizon strategy and Section~\ref{sec:examples} provides illustrative numerical examples. 

\section{System representation}
\label{sec:modelling}
We consider optimal control of a nonlinear system with a discrete time model representation of the form:
\begin{equation}
    g(x_{k+1},x_k,u_k) = 0
\label{eqn:statespace}
\end{equation}
where $x \in \mathbb{R}^{n_x}$ and $u\in \mathbb{R}^{n_u}$ denote the state and control input, and $(x_k,u_k) = (0,0) $ is an equilibrium point. In this section we derive the discrete time system model using a variational approach for both single rate and multirate cases. 


\subsection{Single-rate discrete model formulation}

Consider a mechanical system, whose configuration $q(t)\in \mathbb{R}^{n_q}$ evolves on a configuration manifold $Q$, with associated tangent bundle $TQ$ and a (regular) Lagrangian $\L: TQ \to \mathbb{R}$. The system behaviour can be influenced by control signals $u\in L^{\infty}$ 
which evolve on a control
manifold $U\subset \mathbb{R}^{n_u}$. Thus the motion of the system is governed by the Lagrange-d'Alembert principle~\cite{Reading71}, which enforces the constraint
\begin{equation}
	\label{eq:LdA1}
\delta \int_{t_0}^{t_f} \L(q,\dot{q})\,dt+ \int_{t_0}^{t_f} (f(q,\dot{q},u)\cdot\delta q)\,dt=0
\end{equation}
for all variations $\delta q$ with $\delta q(t_0) = \delta q(t_f) = 0$, where the generalized forces $f(q,\dot{q},u)$ are assumed to be differentiable.

To represent the system dynamics in discrete time we define the discrete path $q_d(t_i) = q_i$ for $t_i = t_0+i\Delta t$, $i =0,...,N$ and approximate the configuration variable using piecewise linear functions : $q_d(t_i) = q_{\,i}+(t-t_i)(q_{\,i+1}-q_{\,i})/\Delta t$, $  \dot{q}(t) \approx (q_{\,i+1}-q_{\,i})/\Delta t $ for $t\in[t_i,t_{i+1}]$. 
Similarly we define a piecewise constant control path $\{u_{i+1/2}\}^{N-1}_{i=0}$ where $ u_{i+1/2} \approx u(t_{i+1/2})$. The action integral on each discrete time step is approximated as 
\[
    \L_{d,i} = \L_{d}(q_i,q_{i+1})\approx\int_{t_i}^{t_{i+1}} \L(q,\dot{q}) \, dt 
\]
 where $\L_d$ is often called the discrete Lagrangian \cite{Reading57}. The virtual work can similarly be approximated as
\[
	f_{i}^{-}\cdot\delta q_{i}+ f_{i}^{+}\cdot\delta q_{i+1} \approx \int^{t_{i+1}}_{t_i}(f(q,\dot{q},u)\cdot\delta q) \, dt .
\]
Here $f_{i}^{-},f_{i}^{+}$ are left and right discrete forces~\cite{Reading71}, with $f_{i-1}^{+}$ ($f_{i}^{-}$) representing the effects of continuous control forces acting over $[t_{i-1},t_{i}]$ ($[t_{i},t_{i+1}]$ respectively) on configuration node $q_i$.
These approximations can be combined to formulate the discrete version of the Lagrange-d'Alembert principle:
\[
    \delta \sum_{i=0}^{N-1} \L_d(q_i,q_{i+1})+ \sum_{i=0}^{N-1} (f_i^- \cdot \delta q_i+f_i^+ \cdot \delta q_{i+1}) = 0
\]
 for all variations $\{\delta q_i\}^N_{i=0}$ vanishing at the endpoints ($\delta q_0=\delta q_N = 0$). The resulting discrete equations of motion  are~\cite{Reading71}
\[
    D_2 \L_d(q_{i-1},q_i) + D_1\L_d(q_{i},q_{i+1})+f_{i-1}^++f_{i}^- = 0
\]
for $i = 1,\ldots,N-1$.
Here $D_n$ denotes the gradient operator with respect to the $n$th argument. These equations define the evolution of the discrete path in configuration space. 

Alternatively we can describe the discrete time system dynamics using an equivalent set of equations 
\begin{subequations}\label{eqn:EOM}
 \begin{align}
     p_i+  D_1\L_d(q_{i},q_{i+1})+ f_{i}^- &= 0
  \label{eqn:EOMa}
\\
     p_{i+1} - D_2\L_d(q_{i},q_{i+1})- f_{i}^+ &=0 
  \label{eqn:EOMb}
 \end{align}
\end{subequations}  
where $\{p_i\}_{i=0}^N$ are the discrete conjugate momenta defined by the discrete forced Legendre transform~\cite{Reading71}. 
%
Defining $x_{i} = [q_i^\top, p_i^\top ]^\top$ and $ n_x = 2n_q$
provides the structure-preserving equations of motion 
in a state space model of the form~(\ref{eqn:statespace}). 

\subsection{Multirate discrete variational model}
\label{sec:multirate_vi}
\begin{figure}[thpb]
      \centering
      \includegraphics[scale=0.24]{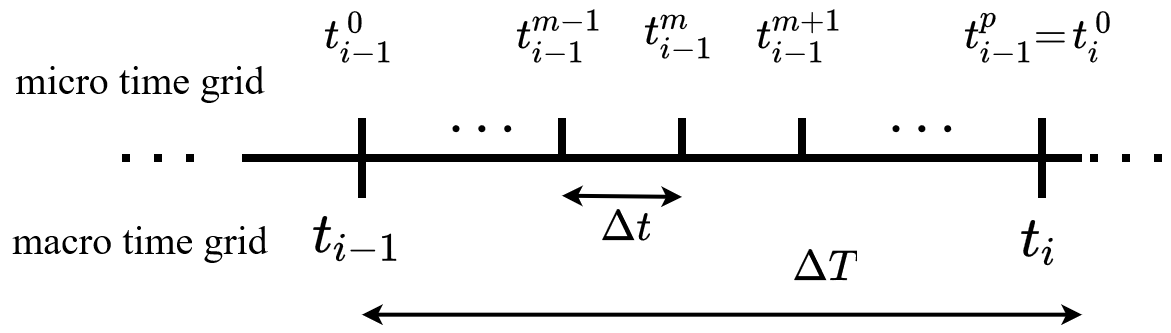}
      \vspace{-2mm}
      \caption{Micro and macro time grid schematic}
      \label{fig:time_grid}
\end{figure}

To derive a multirate discrete time model using the variational approach,
let us assume that the configuration variables can be separated into slow ($q^s\in\mathbb{R}^{n_q^s}$) and fast ($q^f\in\mathbb{R}^{n_q^f}$) components as $y(q) = (q^s, q^f)$, $n_q=n_q^s+n_q^f$, and that the system is under the influence of slow ($u^s\in\mathbb{R}^{n_u^s}$) and fast ($u^f\in\mathbb{R}^{n_u^f}$) controls with $u=(u^s,u^f)$. 
%
We introduce two time grids (Fig.~\ref{fig:time_grid}): a macro time grid with discrete time step $\Delta T$, and a micro time grid obtained by subdividing each macro time step into $\p$ equal steps ($\Delta T = \p\Delta t $). 

We define slow ($\{q_{i}^{s} \}^N_{i=0}$) and fast ($\{q_i^{f,m}\}^{ \mathfrak{p}}_{m=0}\}^{N-1}_{i=0}$) discrete time variables on the macro and micro time grids respectively, and approximate these using piecewise-linear polynomials analogously to the single-rate case. The slow and fast control paths are similarly approximated using piecewise constant functions on the macro and micro time grids respectively as $\{u_{i+1/2}^{s}\}^{N-1}_{i=0}$, $u_{i+1/2}^{s} \approx u^s(t_{i+1/2})$ and $\{ \{u_i^{f,m+1/2}\}^{ \mathfrak{p}-1}_{m=0}\}^{N-1}_{i=0}$, $u_i^{f,m+1/2} \approx u^f(t_i^{m+1/2}) $. 
The discrete multirate Lagrangian is approximated as
\begin{align*}
& \L_{d,i}^{m}(q_i^{s,m},q_i^{s,m+1},q^{f,m}_i,q^{f,m+1}_i) \approx \int^{t_{i}^{m+1}}_{t_{i}^{m}}
\!\!\!\! \L\bigl(y(q),y(\dot{q})\bigr))  \, dt 
\\
& \L_{d,i}(q_i^s,q_{i+1}^s,q^f_i) = \sum^{ \mathfrak{p}-1}_{m=0} \L_d^{i,m} \approx \int^{t_{i+1}}_{t_{i}}
\! \L\bigl(y(q),y(\dot{q})\bigr))  \, dt 
\end{align*}
where $q_i^{s,m} = \frac{ \mathfrak{p}-m}{ \mathfrak{p}} q_i^s + \frac{m}{ \mathfrak{p}} q_{i+1}^s$ and $q_i^f=\{q_i^{f,m}\}^{ \mathfrak{p}}_{m=0}$.
Similarly the virtual work is expressed as 
\begin{align*}
&\sum^{ \mathfrak{p}-1}_{m=0}(f_i^{f,m,-} \!\!\cdot \delta q_i^{f,m}\!+\!f_i^{f,m,+} \!\!\cdot \delta q_i^{f,m+1}) 
\\
&\;\;\;\; +\!    f_i^{s-}\!\!\cdot\delta q_i^s \!+\! f_k^{s+}\!\!\cdot\delta  q_{i+1}^s  \approx\smash[t]{\int^{t_{i+1}}_{t_{i}}} (f^s\cdot \delta q^s+ f^f\cdot \delta q^f) \, dt
\end{align*}
where 
$f\cdot (q^s,q^{f}) = f^{s}\cdot q_{i}+ f^{f}\cdot q^f$. Using these approximations in the discrete Lagrange-d'Alembert principle gives the discrete multirate Euler Lagrange equations~\cite{Reading8, MyASS2020paper}:
\begin{subequations}\label{eqn:DMEOM}
\begin{align}
 &p_i^{s} + D_{q_i^{s}}\L_{d,i}+f_i^{s,-} =0 
\label{eqn:DMEOM1}\\
 &   p_{i+1}^{s} - D_{q_{i+1}^{s}}\L_{d,i} - f_{i}^{s,+} = 0  
\label{eqn:DMEOM2}\\
 &  p_i^{f,0} + D_{q_i^{f,0}}\L_{d,i}^{0}+f^{f,0,-}_i =0
\label{eqn:DMEOM3}\\
 &  p_{i+1}^{f,0} - D_{q_i^{f, \mathfrak{p}}}\L_{d,i}^{\mathfrak{p}-1}-f^{f, \mathfrak{p}-1,+}_{i} =0
\label{eqn:DMEOM4}\\
 & D_{q_i^{f,m}}\L_{d,i}+f_i^{f,m,-}+f_i^{f,m-1,+} = 0
\label{eqn:DMEOM5}
\end{align}
\end{subequations}
where (\ref{eqn:DMEOM1})- (\ref{eqn:DMEOM4}) apply for $i=0,\ldots , N-1$ and (\ref{eqn:DMEOM5})  for $i=0,\ldots , N-1$, $m=1,\ldots,\p$. Here $p_{d}^{s} = \{p_{i}^{s} \}^N_{i=0}$ and  $p_d^f=\{ \{p_i^{f,m}\}^{ \mathfrak{p}}_{m=0}\}^{N-1}_{i=0}$ represent the slow and fast conjugate momenta (\cite{Reading8, MyASS2020paper}).

Equations (\ref{eqn:DMEOM1}-e) can be expressed in the form of~(\ref{eqn:statespace}) with
\begin{equation} \label{eqn:MDMOCx}
x_i = (q_i^s, q_i^{f,0}, p_i^s, p_i^{f,0}), \quad 
u_i = (u_{i+1/2}^s, u_i^{f}) ,
\end{equation}
$n_x = 2n_q^s+2n^f_q$, $n_u = n_u^s+ \p n_u^f$,  and $u^f_i = \{u^{f,m+1/2}_i\}^{ \p-1}_{m=0}$. 
This model formulation is derived by using  (\ref{eqn:DMEOM5}) to determine an expression for $q_i^{f}$, and then eliminating this variable from (\ref{eqn:DMEOM1})-(\ref{eqn:DMEOM4}). Thus the state vector retains the values of the slow and fast variables only on the macro time grid. 
If the fast configuration variables $q_i^f$ need to be evaluated on the micro time grid, these can be obtained using (\ref{eqn:DMEOM5}) or its equivalent couple, for $i = 0,\ldots,N-1$, $m=1,\ldots, \p-1$,
\begin{equation}
\begin{array}{c}
    p_i^{f,m} + D_{q_i^{f,m}}\L_{d,i}^{m}+f^{f,m,-}_i =0\\
    p_{i}^{f,m+1} - D_{q_i^{f,m+1}}\L_{d,i}^{m}-f^{f,m,+}_{i} =0 .
    \label{eqn:DMEOM8}
    \end{array}
\end{equation}

\section{Linearization}
\label{sec:linearization}
The proposed NMPC approach relies on repeatedly linearizing the nonlinear system model~(\ref{eqn:statespace}).
We consider a linear model approximation of the form 
\begin{align}
    0 &= g(x_{k+1},x_k,u_k ) \nonumber\\
    &= M_k x^\delta_{k+1} + D_k x^\delta_k + J_k u^\delta_k + e(x_{k+1},x_k,u_k )
\label{eqn:approx}
\end{align} 
where $x_k^\delta = x_k-x_k^a$, $u^\delta = u_k - u_k^a$, 
$g(x_{k+1}^a,x_k^a,u_k^a )=0$,
$M_k,D_k\in \mathbb{R}^{n_x\times n_x}$, $J_k\in \mathbb{R}^{n_x\times n_u}$ define the linearisation about $x^a_{k+1}x^a_k,u^a_k$,
and $e(x_{k+1},x_k,u_k)$ is the approximation error.

 \subsection{Single-rate variational linearization and error bounds}
For the single-rate case, a linear system approximation can be derived using Jacobian linearization of (\ref{eqn:statespace}). 
Alternatively, we show in this section that a linear discrete variational representation is obtained using (\ref{eqn:EOM}a,b) with a quadratic approximation of the Lagrangian $\L(q,\dot{q})$ and a linear approximation of the generalized forces $f(q,\dot{q},u)$ of the form 
\begin{align}
\tilde{\L}(q,\dot{q},q^a,\dot{q}^a) &\approx {\L}(q^a,\dot{q}^a)+ \nabla_q  \L\vert_{(q^a,\dot{q}^a)}\delta q +\nabla_{\dot{q}}  \L \vert_{(q^a,\dot{q}^a)} \delta\dot{q}
\nonumber\\
& +  \tfrac{1}{2} \delta  q^{\top} \nabla^2_q  \L\vert_{(q^a,\dot{q}^a)}\delta q
 + \tfrac{1}{2} \delta\dot{q}^{\top} \nabla^2_{\dot{q}}  \L \vert_{(q^a,\dot{q}^a)} \delta\dot{q}
  \label{eqn:Lapprox}\\
\tilde{f}(q,\dot{q},u,a) &\approx {f}(a)+  \nabla_q  f\vert_{a}\delta q 
 + 
 \nabla_{\dot{q}}  f \vert_{a} \delta\dot{q}
 +\nabla_{u} f\vert_{a}\delta u
 \label{eqn:fapprox}
\end{align}
Here $a = (q^a,\dot{q}^a,u^a)$ is the approximation point, $\delta q = q-q^a$, $\delta\dot{q} = \dot{q}-\dot{q}^a$, $\delta u = u-u^a$. To simplify presentation we assume that $\L$ and $f$ are separable: $\L(\dot{q},q)$ =$V(\dot{q})+W(q)$, $f(q,\dot{q},u) = z(q)+y(\dot{q})+r(u)$. 
 
\begin{lem}\label{lem:lin}
If $\tilde{\L}(q,\dot{q},q^a,\dot{q}^a) $ is quadratic  in $q,\dot{q}$ and $\tilde{f}(q,\dot{q},u,a)$ is linear in $q,\dot{q},u$, with discrete approximations 
  \begin{align}
\tilde{\L}_{d,i} &= \tilde{\L}_{d}(q_i,q_{i+1},q^a,\dot{q}^a) \nonumber \\
&= \Delta t \, \tilde{\L}(bq_i+cq_{i+1},\tfrac{q_{i+1}-q_{i}}{\Delta t},q^a,\dot{q}^a) \label{eqn:disc_lag_quad}\\
\tilde{f}_{i}^\pm 
&= \tfrac{\Delta t}{2} \tilde{f}( bq_i+cq_{i+1},\tfrac{q_{i+1}-q_{i}}{\Delta t},u_{i+1/2},a) \label{eqn:disc_force_lin}
\end{align}
for some constants $b,c$ then the discrete model (\ref{eqn:EOM}), with $\L_d, f^\pm_i$ replaced by $\tilde{\L}_d, \tilde{f}^\pm_i$, is linear in $q_i,q_{i+1},p_i,p_{i+1}$ and $u_{i+1/2}$.
\end{lem}

\proof 
The discrete approximation $\tilde{f}_{i}^\pm$ inherits the linearity of $\tilde{f}(q,\dot{q},u,a)$. Similarly, since $\tilde{\L}(q,\dot{q},q^a,\dot{q}^a) $ is quadratic, the discrete approximation $ \tilde{\L}_{d,i}$ is also quadratic in $q_i$ and $q_{i+1}$. Thus the gradients  $D_1 \tilde{\L}_d(q_{i},q_{i+1},q^a,\dot{q}^a)$, $D_2 \tilde{\L}_d(q_{i},q_{i+1},q^a,\dot{q}^a)$ are linear in $q_i$, $q_{i+1}$ leaving all terms in (\ref{eqn:EOM}a,b) linear in $q_i,q_{i+1},p_i,p_{i+1},u_{i+1/2}$. 
\endproof

Other appropriate discrete approximations of  $\tilde{\L}$, $\tilde{f}$ can also be used. For example different integration schemes based on the choice of $b,c$ or backward or central differences for the approximation of the velocity.  

 
\begin{thm}\label{thm:lin}
 The linear model obtained at point $a=(q^a,\dot{q}^a,u^a)$ using (\ref{eqn:Lapprox}), (\ref{eqn:fapprox}), (\ref{eqn:disc_lag_quad}), (\ref{eqn:disc_force_lin}) with $\L_d, f_i^\pm$ replaced by $\tilde{\L}_d, \tilde{f}_i^\pm$ in  (\ref{eqn:EOM}) is identical to the Jacobian linearization of (\ref{eqn:EOM}) with $L_{d,i}(q_i,q_{i+1}) = \Delta t L(bq_i+cq_{i+1},\frac{q_{i+1}-q_i}{\Delta t})$, $f^{\pm}_i(q_i,q_{i+1},u_{i+1/2}) = \frac{\Delta t}{2} f(bq_i+cq_{i+1},\frac{q_{i+1}-q_i}{\Delta t},u_{i+1/2})$ around the point $q^a_i,q^a_{i+1},u^a_{i+1/2}$ such that $bq^a_i+cq^a_{i+1}=q^a$, $\frac{q^a_{i+1}-q^a_{i}}{\Delta t}=\dot{q}^a$, $u^a_{i+1/2}=u^a$ for all $i$.
\end{thm}
 
A proof of Theorem~\ref{thm:lin} is given in the Appendix.
The equivalence shown by Theorem~\ref{thm:lin} allows us to obtain important insight into linearized variational models and investigate several of their characteristics, which hold regardless of which of the two linearization methods is used.

\begin{thm}\label{thm:strpr}
If $\L$ and $f$ are nonlinear, then the use of a quadratic Lagrangian approximation and a linear force approximation of the form (\ref{eqn:Lapprox}),(\ref{eqn:disc_lag_quad}) and (\ref{eqn:fapprox}),(\ref{eqn:disc_force_lin}) in (\ref{eqn:EOM}) results in a structure-preserving linear discrete model if the same quadratic Lagrangian approximation is used at each discrete time step, i.e.~if the approximations (\ref{eqn:Lapprox}) and (\ref{eqn:disc_lag_quad}) are computed at the same point $a$ at all times. 
\end{thm}
 \proof
 In the absence of forcing, the scheme (\ref{eqn:EOM}) is symplectic and guarantees the preservation of the discrete momentum map \cite{Reading57}. In the presence of forcing, the preservation of discrete momentum is likewise guaranteed by the Discrete Noether Theorem with Forcing \cite{Reading71}. These properties are maintained by  (\ref{eqn:EOM}) for any continuous force and regular Lagrangian expression and their appropriate discrete approximations \cite{Reading71,Reading57}. Thus the use of (\ref{eqn:fapprox}),(\ref{eqn:disc_force_lin}) and the same Lagrangian function (\ref{eqn:Lapprox}),(\ref{eqn:disc_lag_quad}) throughout the horizon guarantees these properties will be inherited by the resulting model (\ref{eqn:EOM}) and  its linearity follows from Lemma \ref{lem:lin}. 
 \endproof

\begin{rem}\label{rem:lin}
In the numerical experiments presented in Section~\ref{sec:mpc} it is also observed  that simulations using successive linearization of the variational model computed around a seed trajectory, with a different Lagrangian approximation at each discrete time step, provide better numerical behaviour than standard non-geometric approaches. In such cases Theorem~\ref{thm:lin} can be rewritten to state the equivalence of the successive linearization schemes obtained from the variational model using either Jacobian Linearization of (\ref{eqn:EOM}) around the points $(q^a_i,q_{i+1}^a,u^a_{i})$ for $i=0,...N$ or the method proposed by Theorem \ref{lem:lin} at $a_i = (\bar{q}_i^a,\bar{\dot{q}}_i^a,u^a_i )$ for $i=0,...N$ where $\{q^a_i\}^N_{i=0}$, $\{u^a_i\}^{N-1}_{i=0}$ denotes the chosen discrete seed trajectory and $\bar{q}_i^a = bq^a_i+cq^a_{i+1}$,  $\bar{\dot{q}}_i^a = \frac{q_{i+1}^a-q^a_i}{\Delta t }$.
\end{rem}
 

\begin{prop}
Lemma \ref{lem:lin} and Theorems \ref{thm:lin} and \ref{thm:strpr} hold for non-separable expressions for $\L$ and $f$. This results in mixed derivative terms in the approximations (\ref{eqn:Lapprox}), (\ref{eqn:fapprox}), and the proof of Theorem 1 is then analogous to the separable case.
\end{prop}
 
In summary, regardless of whether the linearization method of Lemma~\ref{lem:lin} or Jacobian linearization is used, the resulting variational model (\ref{eqn:EOM}) retains advantageous qualitative properties in its linearized form. Moreover,  from the equivalence of the two linearization approaches, the error $e$ in the Jacobian linearization of (\ref{eqn:EOM}) can be related directly to the error in the approximation of the Lagrangian $e_{\mathcal{\tilde{L}}}$ and the generalized force $e_{\tilde{f}}$ (using the notation of Appendix A) as
 \begin{multline*}
e(q_i,q_{i+1},u_{i+1/2},a) = 
\\
\Delta t \begin{bmatrix}
           \nabla_{q_i}  e_{\tilde{\L}}(\bar{q}_i,\bar{\dot{q}}_i,a) + \frac{1}{2}e_{\tilde{f}^{-}}(\bar{\bar{q}}_i ,\bar{\dot{q}}_i ,u_{i+1/2},a)  \\
                     \nabla_{q_{i+1}}  e_{\tilde{\L}}(\bar{q}_i,\bar{\dot{q}}_i,a) - \frac{1}{2}e_{\tilde{f}^{+}}(\bar{\bar{q}}_i ,\bar{\dot{q}}_i ,u_{i+1/2},a) 
     \end{bmatrix}
 \end{multline*}
%
%
Thus 
bounds on  $e_{\tilde{\L}}$ and $e_{\tilde{f}}$ allow the derivation of bounds for  $e$. 
Hence the proposed linearization approach provides a method of deriving the error bounds required for the tube NMPC approach of Section~\ref{sec:mpc}.
%
The approach is simplified if either the Lagrangian is quadratic and only the force needs to be approximated or the force is linear and only the Lagrangian needs to be approximated. 

\subsection{Multirate variational linearization and error bounds}
\label{sec:multirate_linearization}
Analogously to the single-rate case, for the multirate approach the linearization of the model can again be achieved using two different techniques and bounds on the error of the approximation of (\ref{eqn:DMEOM})-(\ref{eqn:DMEOM8}) can be obtained based on bounds on the approximation of $\L$ and $f$ made at each macro time step.
Linearization of equation (\ref{eqn:DMEOM5}) or its equivalent (\ref{eqn:DMEOM8}) provides expressions for $q_i^f$ and $p_i^f$ that are linear in $x_i$ and $u_i$. Thus, although the proposed state~(\ref{eqn:MDMOCx}) does not contain discrete instances of the configuration and momenta at all micro time nodes, penalties on these variables can be included via a quadratic cost function. Similarly, the linear expressions for $q_i^f$ and $p_i^f$ allow linear constraints to be imposed on the configuration and momentum variables at all micro and macro time steps using sets of linear constraints.

\section{Receding horizon control strategy}
\label{sec:mpc}
This section describes the control problem and proposed NMPC law. We consider optimal regulation of the system~(\ref{eqn:statespace}) subject to a quadratic cost and linear constraints of the form
\begin{align}
&\sum^{\infty}_{k=0}\Big(  \Vert x_k \Vert^2_{Q}+ \Vert u_k \Vert^2_{R}\Big)
\label{eqn:cost}
\\
&F x_k+G u_k \leq h, \quad \forall k \geq 0
\label{eqn:linconstr}
\end{align}
where $ \lVert x \rVert^2_{Q} = x^\top Q x$, $F \in \mathbb{R}^{n_c \times n_x}$ and $G \in \mathbb{R}^{n_c \times n_u}$. A measurement of the state $x_k$ is assumed to be available at each discrete time step $k$, $R\succ 0$ is a positive definite matrix and ${Q \succeq 0}$.
The NMPC strategy is based on~\cite{OrigPaper}, which linearizes the system model about predicted trajectories and uses tubes to bound the effects of approximation errors.

We consider two predicted state and input trajectories of the system~(\ref{eqn:statespace}): $(\mathbf{x}^a_k,\mathbf{u}^a_k)= \{x^a_{i|k},u^a_{i|k}\}_{i = 0}^{N}$ (the `seed' trajectory), and  $(\mathbf{x}_k,\mathbf{u}_k)=\{x_{i|k},u_{i|k}\}_{i \geq 0}$, with $x^a_{0|k}=x_{0|k}=x_k$.
Here $(\cdot)_{i|k}$ indicates the prediction at time $k$ of the value of $(\cdot)$ at time $k+i$ and $N$ is a prediction horizon.
The seed trajectory $(\mathbf{x}^a_k,\mathbf{u}^a_k)$ is exactly known but $(\mathbf{x}_k,\mathbf{u}_k)$ lies in a sequence of sets (i.e.~a tube) determined by the linearization of~(\ref{eqn:statespace}) around $(\mathbf{x}^a_k,\mathbf{u}^a_k)$ and its associated linearization error bounds.


Let the difference between these two trajectories be
\begin{equation}
    x^{\delta}_{i|k} = x_{i|k}-x^a_{i|k} ,
\quad
u^{\delta}_{i|k} = u_{i|k} -u^a_{i|k} ,
    \label{eqn:one}
\end{equation}
then, from (\ref{eqn:approx}), the linearization of (\ref{eqn:statespace}) around  $(\mathbf{x}^a_k,\mathbf{u}^a_k)$ gives
\[
M_{i|k}x^{\delta}_{i+1|k}+ D_{i|k} x^{\delta}_{i|k}+ J_{i|k} u^{\delta}_{i|k}
+ e(x_{i+1|k},x_{i|k},u_{i|k}) = 0
\]
which (assuming that $M_{i|k}$ is invertible) can be expressed as
\begin{equation}\label{eqn:ltv}
    x^{\delta}_{i+1|k}= A_{i|k} x^{\delta}_{i|k}+ B_{i|k} u^{\delta}_{i|k} + w_{i|k}, 
\ \ x^\delta_{0|k} = 0
\end{equation}
where $M_{i|k}A_{i|k}=-D_{i|k}$, $M_{i|k}B_{i|k}=-J_{i|k}$ and where $w_{i|k}$ accounts for linearization error.
We parameterize $u^{\delta}$ in terms of a feedback law and a feedforward control sequence $\{\nu_{i|k}\}_{i=0}^{N-1}$, where each $\nu_{i|k}$ is an optimization variable, as
\begin{equation}\label{eqn:two}
\begin{aligned}
u^{\delta}_{i|k} &= K_{i|k}x_{i|k}^{\delta}+\nu_{i|k}
\\
 x^{\delta}_{i+1|k} &= \Phi_{i|k} x^{\delta}_{i|k}+ B_{i|k} \nu_{i|k} + w_{i|k} 
\end{aligned} \ \ \Biggr\} \ \ i=0,\ldots,N-1
\end{equation}
Here $\Phi_{i|k} = A_{i|k}+B_{i|k}K_{i|k}$, and $K_{i|k}$ is chosen to stabilize the time-varying linear system (\ref{eqn:ltv}) in order to improve the feasibility and numerical conditioning of the NMPC scheme (see e.g.~\cite{Reading854}). 
%
For prediction time steps $i\geq N$ a fixed, locally stabilizing, feedback law is assumed for the linearization of~(\ref{eqn:statespace}) around the target point $(x^a,u^a) = (0,0)$,
\begin{equation}\label{eqn:mode21}
\begin{aligned}
u_{i|k} &= \hat{K}x_{i|k} \\
x_{i+1|k} &= \hat{\Phi} x_{i|k}+  \hat{w}_{i|k}
\end{aligned} \ \ \Biggr\} \ \ i \geq N
\end{equation}
where $\hat{\Phi} = \hat{A}+\hat{B}\hat{K}$ and $\hat{w}_{i|k}$ represents linearization error.

Assuming $g$ is differentiable, \cite{OrigPaper} describes convex polytopic sets bounding $\{w_{i|k}\}_{i=0}^{N-1}$ and $\{\hat{w}_{i|k}\}_{i\geq N}$ of the form 
\begin{align}
    w_{i|k} &\in \Co \bigl( \{ C^j_{i|k} x^{\delta}_{i|k} + D^j_{i|k} u^{\delta}_{i|k}\}_{j = 1}^{n_j} \bigr)
\label{eqn:wbound} \\
    \hat{w}_{i|k} &\in \Co \bigl( \{\hat{C}^j x_{i|k} + \hat{D}^j u_{i|k}\}_{j = 1}^{n_j} \bigr)
\label{eqn:whatbound}
\end{align}
where the matrices $\{\hat{C}^j,\hat{D}^j\}_{j=1}^{n_j}$ are fixed and the matrices $\{\{C^j_{i|k},D^j_{i|k}\}_{j=1}^{n_j}\}_{i=0}^{N-1}$ depend on $\{ {\bf x}_k^a, {\bf u}_k^a \}$.
%
For convenience we decompose $x^\delta$ into nominal and disturbed components,
\[
x^{\delta}_{i|k} = z_{i|k}+\e_{i|k}
\ \ \Biggr\{ \ \ 
\begin{aligned}
z_{i+1|k} &= \Phi_{i|k}z_{i|k}+B_{i|k} \nu_{i|k}, & z_{0|k}=0
\\
\e_{i+1|k} &= \Phi_{i|k}\e_{i|k}+w_{i|k}, & \e_{0|k}=0
\end{aligned}
\]
Then the linearization error bounds (\ref{eqn:wbound}), (\ref{eqn:whatbound}) allow a tube to be constructed containing the uncertain sequence $\{\e_{i|k}\}_{i=0}^{N-1}$. The tube cross-sections are chosen to be ellipsoidal:%
\begin{subequations}\label{eqn:xnconstr}%
\begin{alignat}{2}
    \e_{i|k} &\in \E(V_{i|k}, \beta^2_{i|k}), &\qquad  i &= 1,\ldots,N
\\ 
    x_{i|k} & \in  \E(\hat{V},1), &\qquad  i &\geq N ,
\end{alignat}
\end{subequations}
where $\E(X,\eta^2) =\{y: y^\top Xy\geq \eta^2 \}$, ${X\succ 0}$, is an ellipsoid, $\{V_{i|k}\}_{i=0}^{N-1}$ and $\{\beta_{i|k}\}_{i=0}^{N-1}$ are optimization variables, and $\hat{V}$ is chosen so that $\E(\hat{V},1)$ is positively invariant under~(\ref{eqn:mode21}).

The preceding representation of predicted state and input trajectories enables a dual mode MPC strategy \cite{smpc} to determine optimal perturbations $({\bf x}^\delta_k,{\bf u}_k^\delta)$ on a feasible seed trajectory $({\bf x}_k^a,{\bf u}_k^a)$, then update the seed trajectory via $({\bf x}_k^a,{\bf u}_k^a)\gets ({\bf x}_k^a,{\bf u}_k^a) + ({\bf x}^\delta_k,{\bf u}_k^\delta)$ and repeat this process iteratively until convergence as outlined in \cite{OrigPaper}. 
To define a computationally tractable RHOCP, ensure convergence of the iteration and closed loop stability, a min-max problem is solved at each iteration: 
\begin{equation}
\!\min_{\substack{ \nu_{i|k}, \beta_{i|k} \\ i = 0,\ldots,N-1}} \ \max_{\substack{x_{N|k}\in\E(\hat{V},1)\\\e_{i|k}\in\E(V_{i|k},\beta_{i|k}^2) \\ i = 0,\ldots,N-1}} 
\!\!
\Vert x_{N|k} \Vert^2_{P} + \sum^{N-1}_{i=0} (  \Vert x_{i|k} \Vert^2_{Q} + \Vert u_{i|k} \Vert^2_{R})
    \label{eqn:cost2}
\end{equation}
subject to (\ref{eqn:one}), (\ref{eqn:two}), (\ref{eqn:whatbound}) and additional constraints ensuring (\ref{eqn:linconstr}) and the recursive membership condition $\e_{i|k} \in \E (V_{i|k}, \beta^2_{i|k})$ for $i = 1,\ldots, N$.
This problem can be formulated as a (convex) second order cone program (SOCP) as outlined in \cite{OrigPaper}.

A semidefinite program (SDP) can be formulated and solved offline to obtain the terminal feedback gain $\hat{K}$, and corresponding $\hat{V}$ such that the volume of $\mathcal{E}(\hat{V},1)$ is maximized subject to a bound on of the form $\max_{x\in\mathcal{E}(\hat{V},1)} \Vert x \Vert^2_P \leq \rho$~\cite{smpc}. 
Similarly the terminal weight $P$ can be computed offline by minimizing $\tr (P)$ subject to, for $j = 1,\ldots,n_j$,
\[
P - (\hat{\Phi}+\hat{C}^j+\hat{D}^j\hat{K})^\top P (\hat{\Phi}+\hat{C}^j+\hat{D}^j\hat{K})
\succeq
Q+\hat{K}^\top R \hat{K}
\]
requiring the solution of a SDP.
The sequences $\{K_{i|k}\}_{i=0}^{N-1}$, $\{V_{i|k}\}_{i=0}^{N-1}$ can be computed online by solving a SDP for each $i = N-1,\ldots,0$ (or set equal to $\hat{K}$ and $\hat{V}$ to reduce the computation burden).
 

Detailed analysis of the scheme and an algorithm for its implementation can be found in \cite{OrigPaper}. In this paper we focus on the choice and accuracy of the model (\ref{eqn:ltv}) and propose the use of a variational integrator. As we will demonstrate  in the NMPC context this representation allows for more accurate long-term simulations, presenting new possibilities for error bound computation and most importantly allowing for a multirate formulation which reduces the computational cost of the overall MPC approach.


\section{Multirate receding horizon control}
\label{sec:multirate_mpc}
Using the multirate model of Section~\ref{sec:multirate_vi} we propose a novel receding horizon strategy executed on the macro time grid. Namely, linearization of the model and calculation of finite horizon prediction are performed only at the macro time steps with the RHOCP formulated based on multirate DMOC using a linearized variational model.
This approach allows for an increase in the prediction horizon time step by a factor $ \p$ for negligible sacrifices in accuracy in the resolution of the fast dynamics while allowing for fast control changes 
\cite{MyASS2020paper}. Furthermore for a given horizon length $N$ in the RHOCP, the number of degrees of freedom in control variables reduces from $n_u^s+n_u^f$ to $n_u^s/\p+n_u^f$. Similar reductions are seen in the number of state constraints and the sparse structure of their Jacobian is maintained. Overall this provides significant computational savings compared to a single-rate approach as we demonstrate in Section~\ref{sec:examples}. Additional advantages of the variational structure-preserving representation are the improved qualitative system representation and the option to develop linearization error bounds on the equations of motion using quadratic and linear approximations of the Lagrangian and generalized forces respectively.

\begin{rem} 
The choice of multirate state space representation in (\ref{eqn:statespace}) is not unique. For the specific choice (\ref{eqn:MDMOCx}), the error in the approximation of (\ref{eqn:DMEOM1})-(\ref{eqn:DMEOM4}) depends on $\{\delta q_i^{f,m}\}_{m=1}^{p-1}$ which are not included in the state perturbation vector $x^\delta_i$. Thus, the bounds for this error cannot straightforwardly be expressed in the form of (\ref{eqn:wbound}) or (\ref{eqn:whatbound}) but rather in the form $w_{i|k}\in\Co(\{W^j_{i|k} x^\delta_{i|k} + O_{i|k}^j u^\delta_{i|k} +\sum_{m=1}^{\p-1}X_{i|k}^{j,m}\delta q_{i|k}^{f,m} \}_{j=1}^{n_j})$. However, a bound on $\delta q_{i|k}^{f,m}$ can be obtained in the form $\delta q_{i|k}^{f,m}\in\Co(\{Y_{i|k}^j x^\delta_{i|k} + Z_{i|k}^j u^\delta_{i|k}\}_{j=1}^{n_j}\})$ using (\ref{eqn:DMEOM5}). Combinations of each vertex of this set with each vertex of the set $\{W_{i|k}^j \delta x^M_{i|k}+ O_{i|k}^j \delta u^M_{i|k} \}_{j=1}^{n_j}$ therefore provide a new bounding set definition for $w$ and $\hat{w}$ of the desired form (\ref{eqn:wbound}) or
(\ref{eqn:whatbound}), but with the number of vertices increased from $n_j$ to $n_j^2$.  On the other hand, including $q_{i-1}^f$ in the definition of the model state increases the state dimension $n_x$ from $2(n_q^s+n_q^f)$ to $2n_q^s+( \p+1)n_q^f$, leading again to an increased number of constraints in both offline computations and the online RHOCP. Thus the preferred definition of the discrete state variable $x_i$ is strongly problem-dependent.
\end{rem}

\section{Examples}
\label{sec:examples}
The results in this section are obtained using Python, with the scipy.optimize trust-constr method solver for the NLP formulations and the CVXPY library for the solution of the SDP and SOCP problems \cite{scipy, cvxpy}. In all NLP, SDP, SOCP problems the equations of motion were implemented using automatic differentiation with jax \cite{jax}, rather than analytically derived and implemented. 

\begin{figure}[!tbp]
  \centering
  \begin{minipage}[b]{120pt}
    \includegraphics[width=\textwidth]{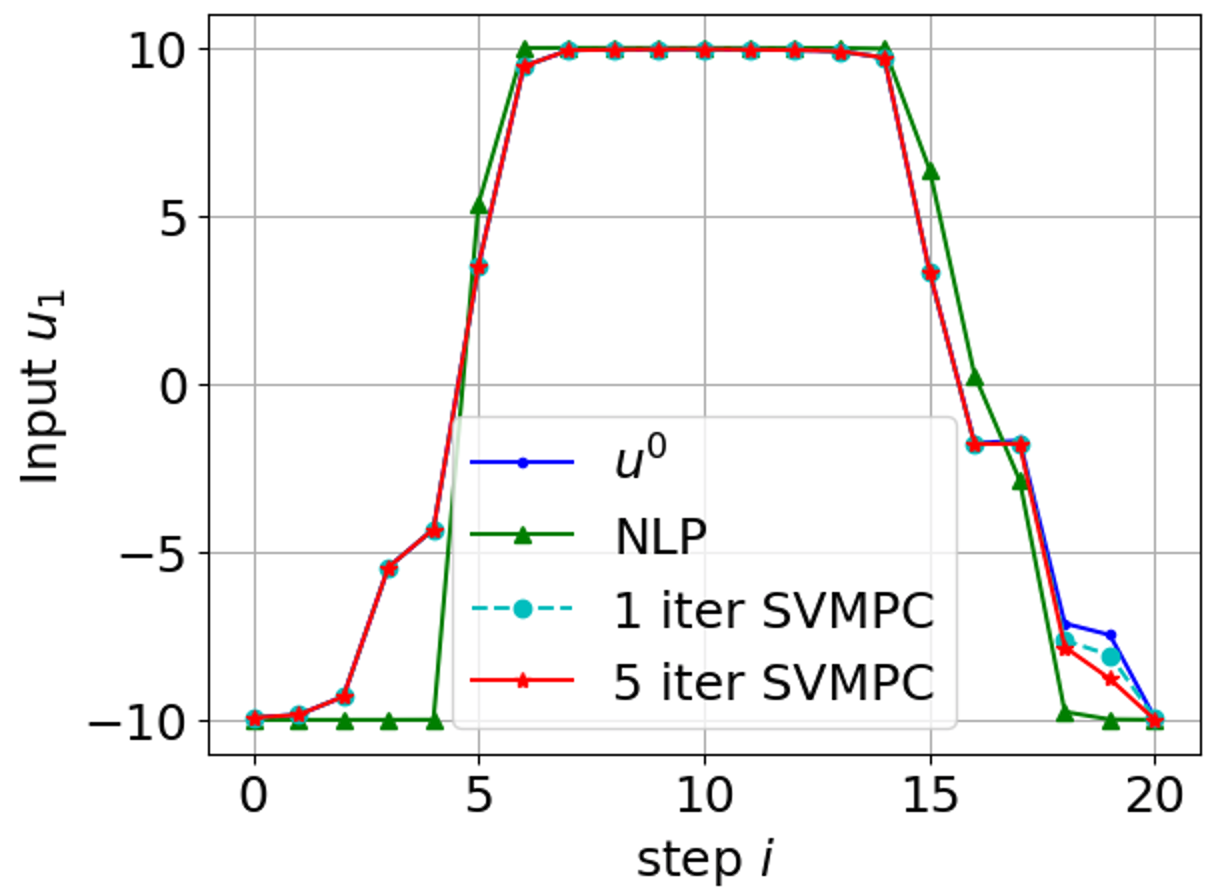}
  \end{minipage}
  \hfil \hspace*{-10pt}
  \begin{minipage}[b]{120pt}
    \includegraphics[width=\textwidth]{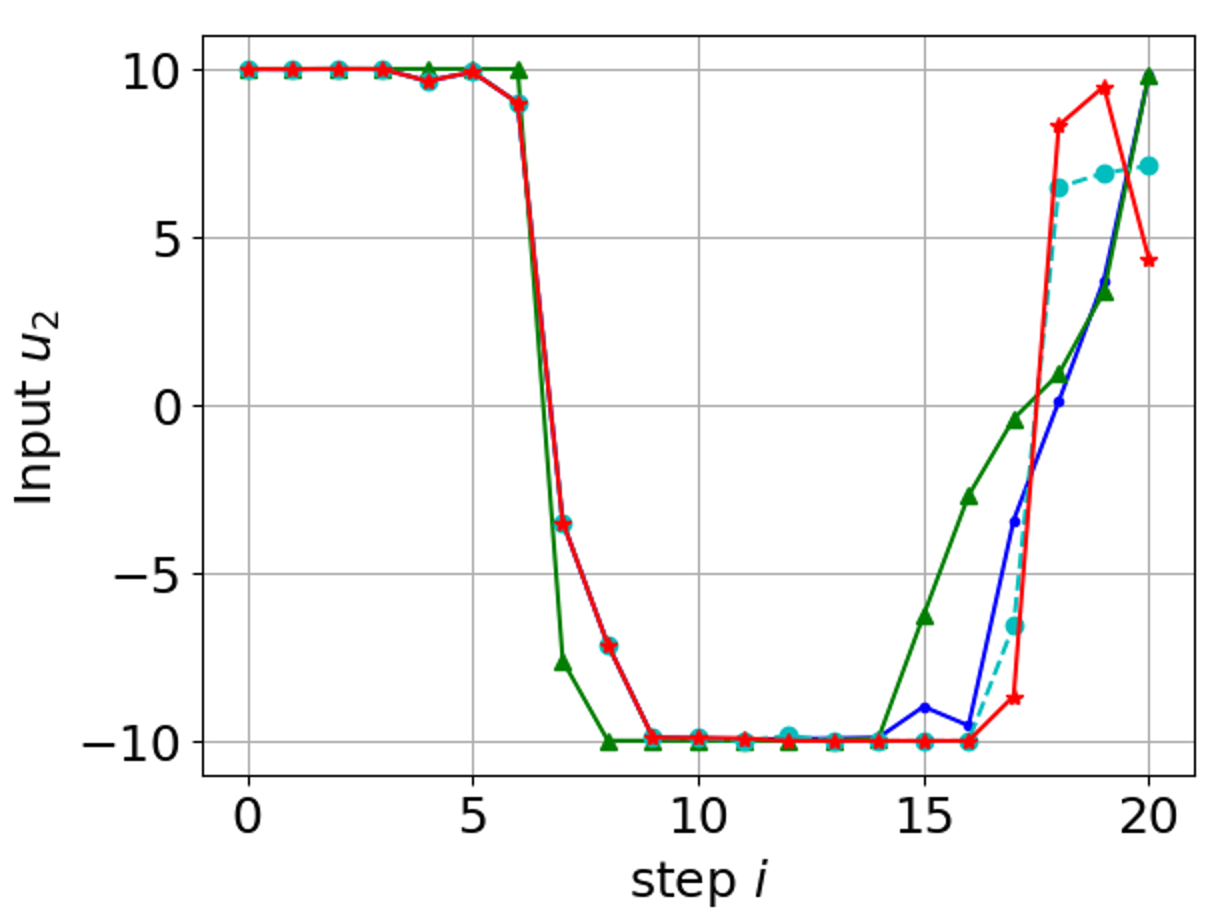}
  \end{minipage}
  \caption{Example control path solution  with $\Delta t = 0.05$, $N=21$}
     \label{fig:single rateSolutions}
\end{figure}

\textit{Example 1}: We first consider a planar model of quadcopter in order to validate the resulting optimal control sequences to those obtained in \cite{OrigPaper}. For this system 
\begin{equation}
\begin{array}{c}
    \L (q,\dot{q}) = \frac{1}{2}\dot{q}^\top \dot{q} - \left[0,g,0\right]^\top \cdot q, \\ f = \left[ (u_1+g)\textrm{sin} \alpha, (u_1+g)\textrm{cos}\alpha, u_2 \right]^\top
    \end{array}
\end{equation}
where $q=\left[y, z, \alpha \right]^\top$ and $u=\left[u_1, u_2\right]^\top$. Here $y,z,\alpha$ represent the horizontal, vertical and angular displacement, $u_1$ and $u_2$ are controls proportional to the net thrust and moment acting on the quadcopter and $g$ is the gravitational acceleration.  To this problem we apply single-rate variational modeling and implement the proposed  MPC scheme with a single-rate variational integrator (SVMPC) using either constant ($V_{i|k}=\hat{V}$) or varying ($V_{i|k}$ computed online) parameters defining MPC tube cross-sections. We refer to these variants as Algorithms 1 and 2, respectively. Fig. \ref{fig:single rateSolutions} presents the predicted control paths for $\Delta t = 0.05$, $N = 21$ and initial condition $x_0 = (0,-1,-1,0,0,0)$. The cost weights are $Q= \diag \{0.1,0.1,10,1,1,1\}$ and $R= \diag \{10^{-4},10^{-3}\}$, and input constraints: $\lvert u_{1,i} \rvert\leq 10$, $\lvert u_{2,i}\rvert\leq 10$. 
The bounds in  (\ref{eqn:whatbound}) are computed using (\ref{eqn:fapprox}) and $\lvert u_{1,i\vert k}\rvert \leq 2.5$, $\lvert u_{2,i \vert k}\rvert \leq 5$, $\lvert \alpha_{i\vert k}\rvert \leq  0.3 $, $\lvert \dot{\alpha}_{i \vert k}\rvert \leq  1$,
the bounds (\ref{eqn:wbound}) are computed online with $\lvert u^\delta_{1,i \vert k} \rvert$,$\lvert u^\delta_{2,i \vert k} \rvert$, $\lvert \alpha_{i \vert k}^\delta\rvert$ and $\lvert \dot{\alpha}_{i \vert k}^\delta\rvert$  restricted to $15\%$ of these bounds. The terminal cost is restricted with $\rho=10$ as in the setting of the original 
investigation in \cite{OrigPaper}. 

   
\begin{figure}[!tbp]
  \centering
  \begin{minipage}[b]{118pt}
    \includegraphics[width=\textwidth]{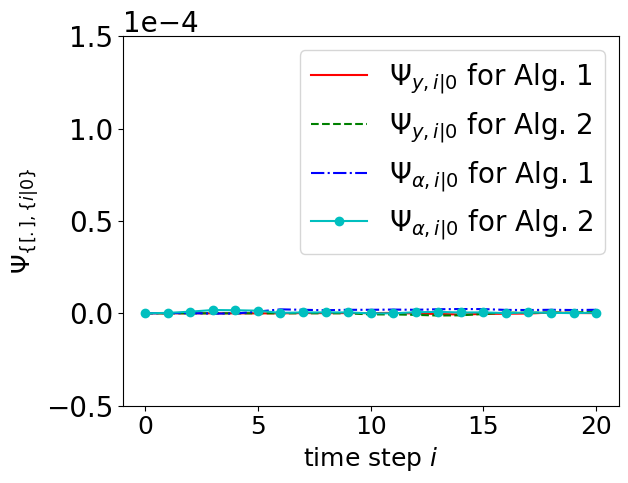}
    \caption{Demonstration of the conservation properties of \textit{Example 1}}
    \label{fig:momentumpreser}
  \end{minipage}
  \hfill
  \begin{minipage}[b]{120pt}
    \includegraphics[width=\textwidth]{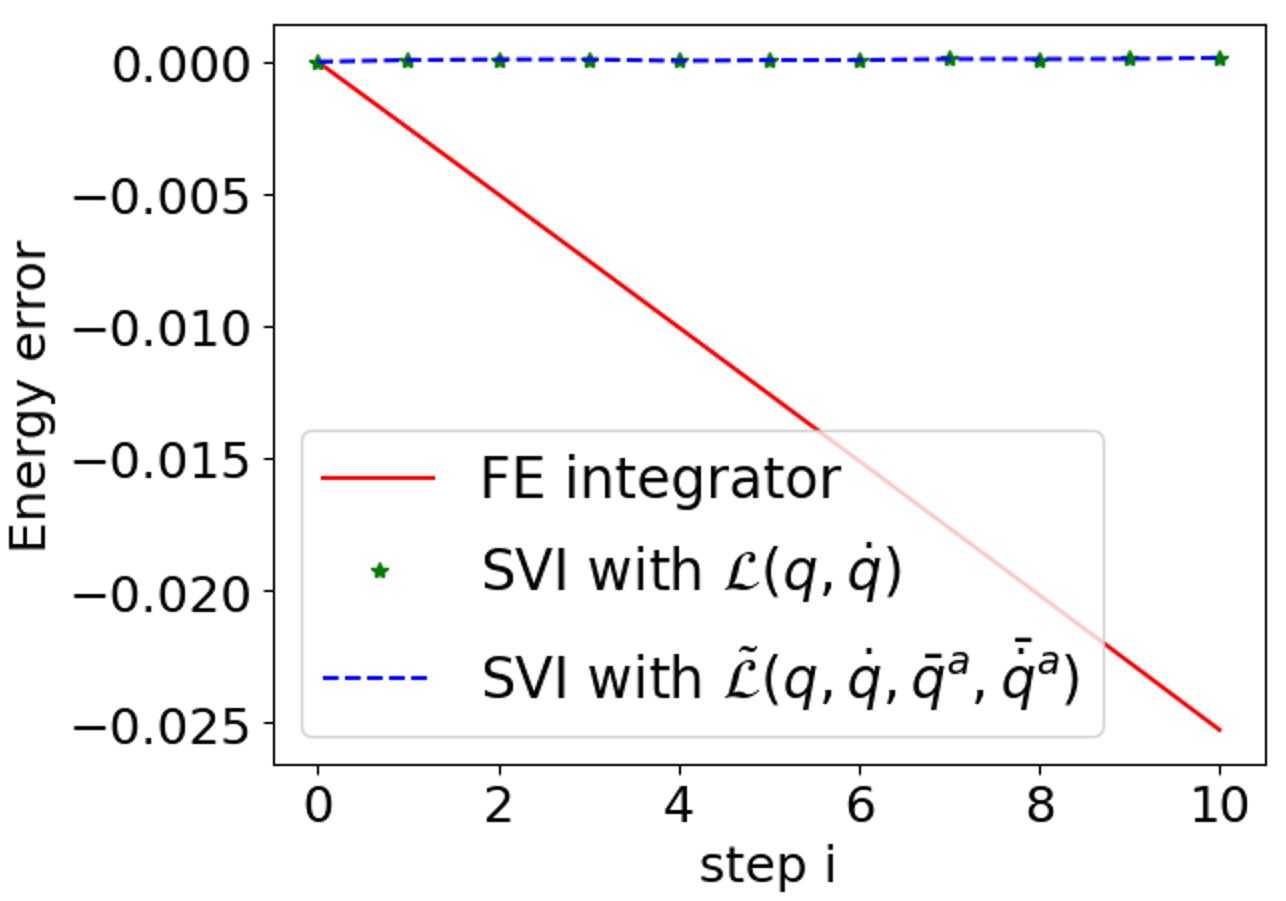}
  \caption{Demonstration of the qualitative properties of the linearized variational scheme for Example 2}
  \label{fig:energypreservation}
  \end{minipage}
\end{figure}

This example also demonstrates the structure-preserving properties of the variational scheme. As the Lagrangian does not explicitly depend on time, $y$ or $\alpha$, the total energy and the momentum associated with $y$ ($p^1)$ and $\alpha$ ($p^{3}$) should remain constant in the absence of external forces  according to Noether’s theorem \cite{Reading57}. In the discrete setting and in the presence of control forces, the momentum map evolves according to the Discrete Noether’s theorem with forcing \cite{Reading71,Reading57}. For \textit{Example 1} this result is equivalent to the following discrete conservation laws:
\begin{align*}
&\Psi_{y,{i\vert k}} = p_{i\vert k}^1- p_{0\vert k}^1-\Delta t \textstyle\sum_{n=0}^{i-1}f^1_{n\vert k}=0\\
&\Psi_{\alpha,{i\vert k}} = p_{i\vert k}^3- p_{0\vert k}^3-\Delta t \textstyle\sum_{n=0}^{i-1}f^3_{n\vert k}=0
\end{align*}
Fig.~\ref{fig:momentumpreser} illustrates these are preserved for the proposed variational approach.

\textit{Example 2:} To investigate the advantages of the multirate approach we introduce the nonlinear Fermi-Pasta-Ulam problem, which has previously been used to investigate multirate schemes \cite{Reading35}. We consider a system consisting of a pair of unit masses connected in series by nonlinear soft and linear stiff springs. 
Let $q^f \in \mathbb{R}^{1}$ represent the length of the stiff spring and $q^s \in \mathbb{R}^{1}$ the location of its center. This choice of coordinates allows for a multirate model with Lagrangian
\begin{equation}
\begin{array}{c}
  \L (q^s,q^f,\dot{q}^s,\dot{q}^f) = \frac{(\dot{q}^s)^2+(\dot{q}^f)^2-(\eta q^f)^2}{2}- \frac{(q^s+q^f)^4+(q^s-q^f)^4}{4} \nonumber 
  \end{array}
\end{equation}
with $f^s = u^s$, $f^f = u^f$ where $u^s$ and $u^f$ represent the resolution of the control applied to the $i$-th mass along $q^s$ and $q^f$ respectively. All tests are performed with $\eta = 50$ and initial conditions  $q^s_0 = \dot{q}_0^s = \dot{q}_0^f =1$ and $q^f_0 = \frac{1}{\eta}$.

To investigate the effect of the proposed linearization scheme on the qualitative properties of the model, simulations were performed in the absence of controls using a Forward Euler (FE) scheme and the single-rate variational scheme (SVI) (\ref{eqn:EOM}) using $\L$ and its approximation $\tilde{\L}$ around a seed trajectory as discussed in Remark~\ref{rem:lin}. 
Results for $t_f = 0.01$ and $\Delta t  = 10^{-3}$ are shown in Fig.~\ref{fig:energypreservation}.
As in \textit{Example 1}, total energy is expected to remain constant in the absence of external forces. The proposed linearization scheme demonstrates more accurate conservation properties than the FE integrator despite the use of a different Lagrangian approximation at each discrete time step $i$ (Fig.~\ref{fig:energypreservation}).

We next consider the proposed multirate variational strategy in an optimal regulation problem starting from the same initial condition with $Q^M= \Delta T \diag\{1,1,1,1\}$, $R^M=  \diag\{\Delta T, \Delta t, \Delta t,\ldots, \Delta t\}$ and constraints $E [q^s, q^f]^\top \leq [10,10]^\top$, $E [\dot{q}^s, \dot{q}^f]^\top \leq [10,10]^\top$, $E [u^s, u^f]^\top\leq [20,20]^\top $ with $E = [1/\sqrt{2}, -1/\sqrt{2};1/\sqrt{2}, 1/\sqrt{2}]$. For this example with $\rho = 1000$, the bounds (\ref{eqn:whatbound}) are based on $E [q^s,q^f]^\top \leq [0.8,0.8]^\top$ and (\ref{eqn:wbound}) are computed using $5\%$ of these bounds.
The computational advantages of the multirate approach are investigated experimentally in a series of tests with a fixed micro time step, $\Delta t=0.01$, and with $t_f=0.6$ and a progressively increasing the micro-macro proportionality $\p$. 
Optimal control paths are obtained using the proposed single-rate and multirate NMPC scheme with approximation of the Lagrangian around a pre-computed seed trajectory for both constant (Algorithm 1) and varying tube cross sections (Algorithm 2) each with \textit{maxiter} $ = 1$. The runtime of the online computation is recorded until discrete time $t_{target}=30\Delta t$. The number of RHOCP steps performed is denoted $n_{rhocp}$. This test is designed to measure the effect of performing the RHOCP more rarely in the multirate approach compared to the single-rate variational scheme. 

\def\strut{\rule{0pt}{6.5pt}}
\begin{table}[]\fontsize{8}{7.2}\selectfont
\caption{\label{table_example}Online computational load }
\vspace{-2mm}
\begin{tabular}{|l|lll|ll|}
\hline\strut
                     & \multicolumn{3}{c|}{Algorithm   1}                                           & \multicolumn{2}{c|}{Algorithm 2}                                             
\\ \hline\strut
$\p $                 & \multicolumn{1}{l|}{1}         & \multicolumn{1}{l|}{3}         & 5         & \multicolumn{1}{l|}{1}         & \multicolumn{1}{l|}{3}                  
\\ \hline\strut
CPU time per         & \multicolumn{1}{l|}{\textbf{3.18}}      & \multicolumn{1}{l|}{\textbf{2.57}}      & \textbf{1.49  }     & \multicolumn{1}{l|}{\textbf{5.52}}       & \multicolumn{1}{l|}{\textbf{4.46} }            \\
iter at   $i=0$      & \multicolumn{1}{l|}{$\pm$0.03} & \multicolumn{1}{l|}{$\pm$0.48} & $\pm$0.02 & \multicolumn{1}{l|}{$\pm$0.17} & \multicolumn{1}{l|}{$\pm$0.01 }   
\\ \hline\strut
$n_{rhocp}$              & \multicolumn{1}{l|}{30}        & \multicolumn{1}{l|}{10}        & 6         & \multicolumn{1}{l|}{30}        & \multicolumn{1}{l|}{10}                
\\ \hline\strut
CPU time to          & \multicolumn{1}{l|}{\textbf{96.3}}     & \multicolumn{1}{l|}{\textbf{27.47}}     & \textbf{11.14}      & \multicolumn{1}{l|}{\textbf{173.79}}     & \multicolumn{1}{l|}{\textbf{51.17} }       \\
reach   $t_{target}$ & \multicolumn{1}{l|}{$\pm$0.24} & \multicolumn{1}{l|}{$\pm$0.43} & $\pm$0.89 & \multicolumn{1}{l|}{$\pm$1.59} & \multicolumn{1}{l|}{$\pm$0.08 }  \\ \hline
\end{tabular}
\end{table}

The results in Table~\ref{table_example} demonstrate that increasing $\p$ reduces the computational load for both approaches. For comparison, the required times for a single NLP solution with $\p =1,3,5$ are $43.57,31.76,29.55$. Thus the proposed single-rate and multirate NMPC schemes outperform the NLP implementation in computation time per discrete time step with the multirate NMPC approach achieving the greatest savings for both Algorithms. For our implementation of Algorithm 1,  the elapsed iteration runtime has a percentage decrease  of $53.14\%$ between $\p= 1$ and $\p= 5$ and the time to reach $t_{target}$ has reduced by a factor of $8.6$. For Algorithm 2, the SDP to compute $V_{i\vert k}$, $K_{i\vert k}$ was not feasible for $\p=5$ for this example. However, $\p=3$ provides percentage decrease of $19.2\%$ per iteration and requires only a third of the time to reach $t_{target}$ compared to the single-rate case.  Of course in general $\p$ cannot be increased indefinitely as this will reduce the resolution of the slow dynamics, and a trade-off therefore exists between computational complexity and accuracy~\cite{MyASS2020paper}. 


\section{CONCLUSIONS}
In this paper we propose a variational discrete time model formulation within a robust tube-based NMPC scheme. Using this model representation we introduce a new method for linearization and error bound extraction. We further propose a novel multirate receding horizon control law which allows the RHOCP optimization to be performed less frequently and allows a reduction in the size of the RHOCP, thereby reducing the number of degrees of freedom, the number of model constraints while retaining a sparse structure for their Jacobian. In numerical examples we show that the proposed multirate variational NMPC scheme can achieve significant runtime savings compared to the single-rate approach for systems exhibiting dynamics and controls on different timescales. We further demonstrate the improved qualitative properties of variational models. 

%



\section*{APPENDIX}
Here we present the proof of Theorem 1. To simplify notation we assume $n_q=n_u=1$. From (\ref{eqn:Lapprox}), $\L$ and $\tilde{\L}$ satisfy
 \begin{equation*}
     \L(q,\dot{q}) =  \tilde{\L}(q,\dot{q},q^a,\dot{q}^a)+ e_{\tilde{\L}}(q,\dot{q},q^a,\dot{q}^a)
 \end{equation*}
\[
     e_{\tilde{\L}}(q,\dot{q},q^a,\dot{q}^a) \!=\!\sum_{n=3}^{\infty}\!\tfrac{1}{n!} \!\Bigl[\!\tfrac{\partial^n \L}{\partial q^n} \vert_{(q^a,\dot{q}^a)} (\delta q)^n +\tfrac{\partial^n \L}{\partial \dot{q}^n} \vert_{(q^a,\dot{q}^a)} (\delta\dot{q})^n\!\Bigr]
\]
based on Taylor series expansion. Using (\ref{eqn:disc_lag_quad}) and
  $   \L_{d,i} = \L_{d}(q_i,q_{i+1})= \Delta t \L(\bar{q}_i,\bar{\dot{q}}_i)$
where $\bar{q}_i = bq_i+cq_{i+1}$, $\bar{\dot{q}}_i = \frac{q_{i+1}-q_{i}}{\Delta t}$, the discrete approximation $\L_{d,i}$ can be defined as 
\begin{align*}     
&\L_{d,i} =  \Delta t \tilde{\L}(\bar{q}_i,\bar{\dot{q}}_i)+ e_{\tilde{\L}_{d,i}}(\bar{q}_i,\bar{\dot{q}}_i,q^a,\dot{q}^a)= \tilde{\L}_{d,i}+ e_{\tilde{\L}_{d,i}}
\\
&e_{\tilde{\L}_{d,i}}=e_{\tilde{\L}_{d,i}}(\bar{q}_i,\bar{\dot{q}}_i,q^a,\dot{q}^a) =\Delta t  e_{\tilde{\L}}(\bar{q}_i,\bar{\dot{q}}_i,q^a,\dot{q}^a)
\\
&\qquad = \sum_{n=3}^{\infty}\tfrac{\Delta t}{n!}\Bigl[\tfrac{\partial^n \L}{\partial \bar{q}_i^n} \vert_{(q^a,\dot{q}^a)} (\delta\bar{q}_i)^n+\tfrac{\partial^n \L}{\partial \bar{\dot{q}}_i^n} \vert_{(q^a,\dot{q}^a)} (\delta\bar{\dot{q}}_i)^n\Bigr].
\end{align*}
Similarly using $f^{\pm}_i= \tfrac{\Delta t}{2} f( \bar{q}_i,\bar{\dot{q}}_i,u_{i+1/2},a)$ and (\ref{eqn:disc_force_lin}) 
\begin{align*}     
&f(q,\dot{q},u) =  \tilde{f}(q,\dot{q},u,a)+ e_{\tilde{f}}(q,\dot{q},u,a), \;\;
f^{\pm}_i =\tilde{f}^{\pm}_i  +e_{\tilde{f}_i^{\pm}} \\
& e_{\tilde{f}_i^{\pm}}= e_{\tilde{f}_i^{\pm}}(\bar{q}_i,\bar{\dot{q}}_i,u_{i+1/2},a) =\tfrac{\Delta t}{2}e_{\tilde{f}}(\bar{q}_i,\bar{\dot{q}}_i,u_{i+\frac{1}{2}},a)
\\
&= \! \sum_{n=2}^{\infty}\! \tfrac{ \Delta t}{2} \! \tfrac{1}{n!}\! \Bigl[\! \tfrac{\partial^n f}{\partial \bar{q}_i^n} \vert_{a} (\delta \bar{q}_i)^n\!\!+\!\tfrac{\partial^n f}{\partial \bar{\dot{q}}_i^n} \vert_{a} (\delta\bar{\dot{q}}_i)^n\!\!+\!\!\tfrac{\partial^n f}{\partial u_{i+1/2}^n} \vert_{a} (\delta u_{i+\frac{1}{2}})^n\!\Bigr]
\end{align*}
Based on the definitions we can rewrite  (\ref{eqn:EOM}a,b) as 
\begin{align*}
     p_i+  D_1(\tilde{\L}_{d,i}+ e_{\tilde{\L}_{d,i}})+ \tilde{f}_{i}^-+e_{\tilde{f}_i^-} &= 0 
\\
     p_{i+1} - D_2(\tilde{\L}_{d,i}+ e_{\tilde{\L}_{d,i}})- f_{i}^+ -e_{\tilde{f}_i^+} &=0 .
 \end{align*} 
  Thus the approximation error $\tilde{e}$ when using using the approximation (\ref{eqn:Lapprox}), (\ref{eqn:fapprox}) in (\ref{eqn:EOM}a,b) can be expressed as
\[
\tilde{e}(q_i,q_{i+1},u_{i+1/2},a) =  \begin{bmatrix}
           D_{q_i} e_{\tilde{\L}_{d,i}}+ e_{\tilde{f}^{-}}  \\
                     -D_{q_{i+1}} e_{\tilde{\L}_{d,i}} - e_{\tilde{f}^{+}} 
     \end{bmatrix} .
 \]
 On the other hand the error of the Jacobian linearization of (\ref{eqn:EOM}) around the point $q^a_i,q^a_{i+1},u^a_{i+1/2}$ such that $u^a_{i+1/2}=u^a$, $bq^a_i+cq^a_{i+1}=q^a$, $\frac{q^a_{i+1}-q^a_{i}}{\Delta t}=\dot{q}^a$ for all $i$ takes the form
\[
    \sum_{j=2}^{\infty}\tfrac{\Delta t}{j!}\Bigl[ b\tfrac{\partial^j}{\partial q^j}\tfrac{\partial \L}{\partial q} \vert_{x^a}(\delta\bar{q}_i)^j- \tfrac{1}{\Delta t}\tfrac{\partial^j}{\partial \dot{q}^j}\tfrac{\partial \L}{\partial \dot{q}} \vert_{x^a} (\delta\bar{q}_i)^j\Big)+e_{\tilde{f}^{-}} 
\]
since $D_1\L_{d,i} =\Delta t \bigl[ b \tfrac{\partial \L}{\partial q} \vert_{(q^a,\dot{q}^a)}-\tfrac{1}{\Delta t} \tfrac{\partial \L}{\partial q} \vert_{(q^a,\dot{q}^a)}\bigr]$
and $x^a=(q^a,\dot{q}^a)$. Comparing this with 
\begin{multline}\label{eqn:label1} 
 D_{q_i} e_{\tilde{\L}_{d,i}}+ e_{\tilde{f}^{-}} 
= \sum_{n=3}^{\infty} \tfrac{\Delta t}{(n-1)!}\Bigl[ b\tfrac{\partial^n \L}{\partial \bar{q}_i^n} \vert_{(q^a,\dot{q}^a)} (\delta\bar{q}_i)^{n-1}\\-\tfrac{1}{\Delta t} \tfrac{\partial^n \L}{\partial \dot{q}^n} \vert_{(q^a,\dot{q}^a)} (\delta\bar{\dot{q}})^{n-1}\Bigr] + e_{\tilde{f}^{-}}
\end{multline} 
it can be seen that the Jacobian linearization is equivalent to the approximation 
 (\ref{eqn:label1}). Similar expansions of (\ref{eqn:EOMb}) and $-D_{q_{i+1}} e_{\tilde{\L}_{d,i}} - e_{\tilde{f}^{+}} $ show that they are also equivalent, thus completing the proof of Theorem 1. 

\bibliographystyle{IEEEtran}
\bibliography{bibliography}

\addtolength{\textheight}{-12cm}   

\end{document}